\documentclass[12pt,leqno]{article}
\usepackage{amsmath, amsthm, amsfonts, amssymb, color}
\usepackage{mathrsfs}
\theoremstyle{definition}

\newcommand{\scr}[1]{\mathscr #1}
\definecolor{wco}{rgb}{0.5,0.2,0.3}

\numberwithin{equation}{section} \theoremstyle{remark}

\newcommand{\ua}{\uparrow}

\title{{\bf Derivative Formula and Applications for Hyperdissipative Stochastic Navier-Stokes/Burgers  Equations}\footnote{Supported in
 part by WIMCS and NNSFC(10721091).}
}
\author{
{\bf Feng-Yu Wang$^{a),b)}$ and Lihu Xu$^{c)}$}\\
\footnotesize{$^{a)}$ School of Math. Sci. and Lab. Math. Com. Sys.,
Beijing Normal
University, Beijing 100875, China}\\
 \footnotesize{$^{b)}$ Department of Mathematics,
Swansea University, Singleton Park, SA2 8PP, UK}\\
\footnotesize{Email: wangfy@bnu.edu.cn;
F.Y.Wang@swansea.ac.uk}\\
\footnotesize{$^{c)}$ PO Box 513, EURANDOM, 5600 MB  Eindhoven. The Netherlands} \\
\footnotesize{Email: xu@eurandom.tue.nl}\\
}
\begin{document}
\def\R{\mathbb R}  \def\ff{\frac} \def\ss{\sqrt} \def\B{\mathbf
B}
\def\N{\mathbb N} \def\kk{\kappa} \def\m{{\bf m}}
\def\dd{\delta} \def\DD{\Delta} \def\vv{\varepsilon} \def\rr{\rho}
\def\<{\langle} \def\>{\rangle} \def\GG{\Gamma} \def\gg{\gamma}
  \def\nn{\nabla} \def\pp{\partial} \def\EE{\scr E}
\def\d{\text{\rm{d}}} \def\bb{\beta} \def\aa{\alpha} \def\D{\scr D}
  \def\si{\sigma} \def\ess{\text{\rm{ess}}}
\def\beg{\begin} \def\beq{\begin{equation}}  \def\F{\scr F}
\def\Ric{\text{\rm{Ric}}} \def\Hess{\text{\rm{Hess}}}
\def\e{\text{\rm{e}}} \def\ua{\underline a} \def\OO{\Omega}  \def\oo{\omega}
 \def\tt{\tilde} \def\Ric{\text{\rm{Ric}}}
\def\cut{\text{\rm{cut}}} \def\P{\mathbb P} \def\ifn{I_n(f^{\bigotimes n})}
\def\C{\mathbb C}      \def\aaa{\mathbf{r}}     \def\r{r}
\def\gap{\text{\rm{gap}}} \def\prr{\pi_{{\bf m},\varrho}}  \def\r{\mathbf r}
\def\Z{\mathbb Z} \def\vrr{\varrho} \def\ll{\lambda}
\def\L{\scr L}\def\Tt{\tt} \def\TT{\tt}\def\II{\mathbb I}
\def\i{{\rm in}}\def\Sect{{\rm Sect}}\def\E{\mathbb E} \def\H{\mathbb H}
\def\M{\scr M}\def\Q{\mathbb Q} \def\texto{\text{o}} \def\LL{\Lambda}
\def\Rank{{\rm Rank}} \def\B{\scr B}
\def\T{\mathbb T}\def\i{{\rm i}} \def\ZZ{\hat\Z}

\maketitle
\begin{abstract} By using  coupling method, a Bismut type derivative  formula is established for the Markov semigroup associated to a class of hyperdissipative stochastic Navier-Stokes/Burgers    equations. As applications, gradient estimates, dimension-free Harnack inequality, strong Feller property, heat kernel estimates and some properties of the invariant probability measure are derived.
\end{abstract} \noindent

 AMS subject Classification:\ 60J75, 60J45.   \\
\noindent
 Keywords: Bismut formula, coupling,  strong Feller, stochastic Navier-Stokes equation.
 \vskip 2cm

\section{Introduction}

Let $H$ be the divergence free sub-space of $L^2(\T^d;\R^d)$, where
$\T^d:=(\R/[0,2\pi])^d$ is the $d$-dimensional torus. The
$d$-dimensional  Navier-Stokes equation (for $d\ge 2$) reads

$$\d X_t= \{\nu \DD X_t-B(X_t,X_t)\}\d t,$$ where $\nu>0$ is the
viscosity constant and $B(u,v):=\mathbf P (u\cdot \nn)v$ for
$\mathbf P: L^2(\T^d;\R^d)\to H$ the orthogonal projection (see e.g.
\cite{T}). When $d=1$ and $H= L^2(\T^d;\R^d)$, this equation reduces
to the Burgers equation. In recent years, the stochastic
Navier-Stokes equations have been investigated intensively, see e.g.
\cite{HM} for the ergodicity of 2D Navier-Stokes equations with
degenerate noise, and see \cite{DD, FR, RX} for the study of 3D
stochastic Navier-Stokes equations. The main purpose of this paper
is to establish the Bismut type derivative formula for the Markov
semigroup associated to stochastic Navier-Stokes type equations, and
as applications, to derive gradient estimates, Harnack inequality,
and strong Feller property for the semigroup.

We shall work with a more general framework as in \cite{LR}, which
will be reduced to a class of  hyperdissipative (i.e. the Laplacian has a power larger than $1$)  stochastic Navier-Stokes/Burgers    equations in Section 2.

Let $(H,\<\cdot,\cdot\>, \|\cdot\|_H)$ be a separable real Hilbert
space, and $(L,\D(L))$ a positively definite self-adjoint operator
on $H$ with $\ll_0:=\inf\si(L)>0,$ where $\si(L)$ is the spectrum of
$L$. Let $V=\D(L^{1/2})$, which is a Banach space  with norm
$\|\cdot\|_V:= \| L^{1/2}\cdot\|.$ Let $Q$ be a Hilbert-Schmidt
linear operator on $H$ with $\text{Ker Q}=\{0\}.$ Then $\D(Q^{-1}):=
Q(H)$ is a Banach space with norm $\|x\|_Q:=\|Q^{-1}x\|_H.$ In
general, for $\theta>0$, let $V_\theta= \D(L^{\theta/2})$ with norm
$\|L^{\theta/2}\cdot\|_H.$ We assume that there exist two constants
$\theta\in (0,1]$ and $K_1>0$ such that $V_\theta\subset \D(Q^{-1})$
and

\ \newline {\bf (A0)}\  $\|u\|_Q^2\le K_1\|u\|_{V_\theta}^2,\ \ u\in
V_\theta.$

\ \newline Moreover, let

$$B: V\times V\to H$$ be a bilinear map such that

\ \newline {\bf (A1)} \  $\<v, B(v,v)\>=0,\ \ v\in V;$

\ \newline {\bf (A2)} \ There exists a constant $C>0$ such that
$\|B(u,v)\|_H^2\le C\|u\|_H^2\|v\|_V^2,\ u,v\in V$;

\ \newline {\bf (A3)} \  There exists a constant $K_2>0$ such that
$\|B(u,v)\|_Q^2\le K_2 \|u\|_{V_\theta}^2\|v\|_{V_\theta}^2,\ \
u,v\in V.$

\

Finally, let $W_t$ be the cylindrical Brownian motion on $H$. We
consider the following stochastic differential equation on $H$:

\beq\label{E} \d X_t= Qd W_t-\{LX_t +B(X_t)\}\d t,\end{equation}
where $B(X_t):=B(X_t,X_t).$ According to \cite{LR}, for any initial value $X_0\in
H$ the equation (\ref{E}) has a unique strong solution, which gives
rise to a Markov process on $H$ (see Appendix for details). For any
$x\in H$, let $X_t^x$ be the solution starting at $x$. Let
$\B_b(H)$ be the set of all bounded measurable functions on $H$. Then

$$P_t f(x):= \E f(X_t^x),\ \ x\in H, t\ge 0, f\in \B_b(H)$$ defines a Markov semigroup $(P_t)_{t\ge 0}.$

We shall adopt a coupling argument to establish a Bismut type derivative formula for  $P_t$, which will imply explicit gradient estimates
and the dimension-free Harnack inequality in the sense of
\cite{W97}. This type of Harnack inequality has been applied to the
study of several models of SDEs and SPDEs, see e.g. \cite{DRW09, L,
LW08, ORW, OY, W07} and references within.

For $f\in \B_b(H), h\in V_\theta, x\in H$ and $t>0$, let

$$D_h P_t f(x)=\lim_{\vv\to 0} \ff 1 \vv \big\{P_t f(x+\vv h)-P_t
f(x)\big\}$$ provided the limit in the right-hand side exists. Let
$\tt B(u,v)= B(u,v)+B(v,u)$.

\beg{thm} \label{T1.1} Assume that {\bf (A0)}-{\bf (A3)} hold for
some constants $\theta\in (0,1], K_1,K_2,C>0.$ Then for any $t>0,
h\in V_\theta$ and $f\in \B_b(H)$, $D_h P_tf$ exists on $H$ and
satisfies

\beq\label{B} D_h P_tf(x)= \E \bigg\{f(X_t^x) \int_0^t \Big\<Q^{-1}
\Big(\ff 1 t \e^{-sL} h- \ff{t-s}t \tt B(X_s^x, \e^{-sL}h)\Big), \d
W_s\Big\>\bigg\},\ \ x\in H.\end{equation}\end{thm}

Let $V_\theta^*$ be the dual space of $V_\theta$. According to Theorem \ref{T1.1}, under assumptions
{\bf (A0)}-{\bf (A3)} we may define the gradient $DP_tf: H\to V_\theta^*$ by letting

$$_{V_\theta^*}\<DP_tf(x), h\>_{V_\theta} = D_h P_t f(x),\ \ x\in H, h\in V_\theta.$$ We shall estimate

$$\|D P_t f(x)\|_{V_\theta^*}:= \sup_{\|h\|_{V_\theta}\le 1} |D_h
P_t f(x)|,\ \ \ x\in H.$$
To this end, let $\|Q\|$ and $\|Q\|_{HS}$ be the operator norm and the
Hilbert-Schmidt norm of $Q: H\to H$ respectively.

\beg{cor}\label{C1.2} Under assumptions of Theorem \ref{T1.1}.
\beg{enumerate} \item[$(1)$] For any $t>0,x\in H$ and $f\in
\B_b(H)$,

$$\|DP_t f(x)\|_{V_\theta^*}^2 \le (P_tf^2(x)) \Big\{\ff{2K_1}t
+\ff{4K_2}{\ll_0^{2-\theta}}\big(\|x\|_H^2
+\|Q\|_{HS}^2t\big)\Big\}.$$\item[$(2)$] Let $f\in \B_b(H)$ be
positive. For any $x\in H, t>0$ and  $\dd\ge 4\ss{K_2}\,
\|Q\|\ll_0^{(\theta-3)/2},$ \beg{equation*}\beg{split}
\|DP_tf(x)\|_{V_\theta^*}\le
&\dd\big\{P_t(f\log f)-(P_tf)\log P_tf\big\}(x)\\
&+\ff 2 \dd \Big\{\ff{K_1}t
+\ff{2K_2}{\ll_0^{1-\theta}}\big(\|x\|_H^2
+\|Q\|_{HS}^2t\big)\Big\}P_tf(x).\end{split}\end{equation*}
\item[$(3)$] Let $\aa>1, t>0$ and $f\ge 0$. The Harnack inequality
$$(P_tf(x))^\aa\le
(P_tf^\aa(y))\exp\bigg[\ff{2\aa\|x-y\|_{V_\theta}^2}{\aa-1}
\Big\{\ff{K_1}t +\ff{2K_2}{\ll_0^{1-\theta}}
\big(\|x\|_H^2\lor\|y\|_H^2+ \|Q\|_{HS}^2t\big)\Big\}\bigg]$$ holds
for $x,y\in H$ such that
$$\|x-y\|_{V_\theta}\le
\ff{(\aa-1)\ll_0^{(3-\theta)/2}}{4\aa\|Q\|\ss{K_2}}.$$ In
particular, $P_t$ is $V_\theta$-strong Feller, i.e.
$\lim_{\|y-x\|_{V_\theta}\to 0} P_t f(y)=P_tf(x)$ holds for $f\in\B_b(H), t>0, x\in H.$
\end{enumerate} \end{cor}

As applications of the Harnack inequality derived above, we have the
following result.

\beg{cor} \label{C1.3} Under assumptions of Theorem \ref{T1.1}.
$P_t$ has an  invariant probability measure $\mu$ such that
$\mu(V)=1$ and hence, $\mu(V_\theta)=1.$ If moreover $\theta\in
(0,1)$, then: \beg{enumerate}
\item[$(1)$] $P_t$ has a unique invariant probability measure $\mu$,
and the measure has full support on $V_\theta$.
\item[$(2)$] $P_t$ has a density $p_t(x,y)$ on $V_\theta$ w.r.t.
$\mu$. Moreover, let
$r_0=\ff{(\aa-1)\ll_0^{(3-\theta)/2}}{4\aa\|Q\|\ss{K_2}}$ and
$B_\theta(x,r_0)= \{y: \|y-x\|_{V_\theta}\le r_0\}$,

\beg{equation*}\beg{split} &\bigg(\int_{V_\theta}
p_t(x,y)^{(\aa+1)/\aa}\mu(\d y)\bigg)^{\aa}\\
& \le \ff 1 {\int_{B_\theta(x,r_0)}
\exp\big[-\ff{2\aa\|x-y\|_{V_\theta}^2}{\aa-1}\big\{\ff{K_1}t
+\ff{2K_2}{\ll_0^{1-\theta}} (\|x\|_H^2\lor\|y\|_H^2+
\|Q\|_{HS}^2t)\big\}\big]\mu(\d y)}<\infty\end{split}\end{equation*}
holds for any $t>0,\aa>1$ and $x\in V_\theta$.  \end{enumerate}
\end{cor}

Note  that the Harnack inequality presented in Corollary \ref{C1.2}
is local in the sense that $\|x-y\|_{V_\theta}$ has to be bounded
above by a constant. To derive a global Harnack inequality, we need
to extend the gradient-entropy inequality in Corollary \ref{C1.2}
(2) to all $\dd>0$. In this spirit, we have the following result.

\beg{thm}\label{T1.4} Under assumptions of Theorem \ref{T1.1}.
\beg{enumerate}
\item[$(1)$] For any $\dd>0$ and any positive
$f\in\B_b(H)$,

 \beg{equation*}\beg{split}
\|DP_tf(x)\|_{V_\theta^*}\le
&\dd\big\{P_t(f\log f)-(P_tf)\log P_tf\big\}(x)\\
&+\ff 2 \dd \Big\{\ff{K_1}{t\land t_\dd}
+\ff{2K_2\e}{\ll_0^{1-\theta}}\big(\|x\|_H^2
+\|Q\|_{HS}^2t\big)\Big\}P_tf(x),\ \  x\in H,
t>0\end{split}\end{equation*}holds for $t_\dd:= \ff{\dd^2
\ll_0^{3-\theta}}{4\|Q\|^2\e K_2}.$
\item[$(2)$] Let $\aa>1, t>0$ and $f\ge 0$. Then
\beg{equation*}\beg{split} (P_tf(x))^\aa\le
(P_tf^\aa(y))\exp\bigg[\ff{2\aa\|x-y\|_{V_\theta}^2}{\aa-1} &\Big\{
K_1\Big(\ff 1 t \lor \ff{4\aa^2\|Q\|^2 \e K_2
\|x-y\|_{V_\theta}^2}{(\aa-1)^2\ll_0^{3-\theta}}\Big)\\
&+\ff{2K_2\e}{\ll_0^{1-\theta}} \big(\|x\|_H^2\lor\|y\|_H^2+
\|Q\|_{HS}^2t\big)\Big\}\bigg]\end{split}\end{equation*} holds for
all $x,y\in \H$.
\end{enumerate}\end{thm}

The remainder of the paper is organized as follows. We first
consider in Section 2 a class of stochastic Navier-Stokes type
equations to illustrate our results, then prove these results in
Section 3.

\section{Stochastic hyperdissipative Navier-Stokes/Burgers
equations}

Let $\T^d= (\R/[0,2\pi])^d$ for $d\ge 1$. Let $\DD$ be the Laplace
operator on $\T^d$. To formulate $\DD$ using spectral
representation, we first consider the complex $L^2$ space
$L^2(\T^d;\C^d)$. Recall that for $a=(a_1,\cdots, a_d),
b=(b_1,\cdots, b_d)\in \C^d$, we have $a\cdot b= \sum_{i=1}^d
a_i\bar b_i$. Let

$$e_k(x)= (2\pi)^{-d/2} \e^{\i (k\cdot x)},\ \ k\in\Z^d, x\in \T^d.$$
Then $\{e_k: k\in\Z^d\}$ is an ONB of $L^2(\T^d;\C)$. Obviously, for a sequence $\{u_k\}_{k\in \Z^d}\subset \C^d$,

$$u:=\sum_{k\in \Z^d} u_k e_k\in L^2(\T^d;\R^d)$$ if and only if $\bar u_k= u_{-k}$ holds for any $k\in \Z^d$ and
$\sum_{k\in \Z^d} |u_k|^2<\infty.$ By spectral
representation, we may characterize $(\DD,\D(\DD))$ on $L^2(\T^d;\R^d)$ as
follows: \beg{equation*} \beg{split} &\DD u= -\sum_{k\in\Z^d}|k|^2
u_k
e_k,\ \ u:=\sum_{k\in\Z^d}u_ke_k\in  \D(\DD),\\
&\D(\DD):=\bigg\{\sum_{k\in\Z^d}u_k\e_k:\ u_k\in\C^d, \bar u_k=u_{-k},
\sum_{k\in\Z^d}
|u_k|^2|k|^4<\infty\bigg\}.\end{split}\end{equation*}

To formulate the Navier-Stokes/Burgers type equation, when $d\ge 2$ we consider the sub-space
divergence free elements of $L^2(\T^d;\R^d)$. It is easy to see that a
smooth vector field

$$u=\sum_{k\in\Z^d} u_k e_k$$ is divergence free if and only if
$u_k\cdot k=0$ holds for all $k\in \Z^d.$ Moreover, to make the
spectrum of $-\DD$ strictly positive, we shall not consider non-zero
constant vector fields. Therefore, the Hilbert space we are working
on becomes

$$H:=\bigg\{\sum_{k\in\ZZ^d} u_k e_k:\, u_k\in
\C^d, (d-1)(u_k\cdot k)=0, \bar u_k=u_{-k}, \sum_{k\in \ZZ^d} |u_k|^2<\infty\bigg\},$$ where $\ZZ^d=\Z^d\setminus\{0\}.$
Since when $d=1$ the condition $(d-1)(u_k\cdot k)=0$ is trivial, the divergence free restriction does not apply for the one-dimensional case.

Let $(A,\D(A))=(-\DD,\D(\DD))|_H$, the restriction of
$(\DD,\D(\DD))$ on $H$, and let $\mathbf P: L^2(\T^d; \R^d)\to H$ be
the orthogonal projection. Let $$L=\ll_0A^{\dd+1}$$ for some
constants $\ll_0,\dd>0$. As in Section 1, define $V=\D(L^{1/2})$ and
$V_\theta=\D(L^{\theta/2}).$ Then

$$B: V\times V\to H;\ \ B(u,v)= \mathbf P (u\cdot \nn)v$$ is a
continuous bilinear (see the (b) in the proof of Theorem \ref{T2.1}
below).  Let $Q=A^{-\si}$ for some $\si>0$, and let $W_t$ be the
cylindrical Brownian motion on $H$. Obviously, $\|Q\|\le 1$ and when
$\si>\ff d 4$,

$$\|Q\|_{HS}^2\le \sum_{k\in\ZZ^d} |k|^{-4\si}<\infty.$$We consider the stochastic differential equation

\beq\label{NS} \d X_t= Q \d W_t -(LX_t+B(X_t))\d t,\end{equation}
where $B(u):=B(u,u)$ for $u\in V$. Thus, we are working on the stochastic hyperdissipative Navier-Stokes (for $d\ge 2$) and Burgers (for $d=1$) equations.

\beg{thm}\label{T2.1} Let $\dd>\ff d 2, \si\in (\ff d 4, \ff \dd 2]$ and
$\theta\in [\ff{2\si+1}{\dd+1},1].$ Then all assertions in Section 1 hold for $K_1= \ff 1 {\ll_0^\theta}$ and

$$K_2= \ff {4^{2\dd\theta +1}} {\ll_0^{2\theta}} \sum_{k\in\ZZ^d}
|k|^{-2(\dd+1)\theta}  <\infty.$$\end{thm}

\beg{proof} Since $\si>\ff d 4$, $Q:H\to H$ is Hilbert-Schmidt. By
Theorem \ref{T1.1} and its consequences, it suffices to verify
assumptions {\bf (A0)}-{\bf(A3)}. Since {\bf (A1)} is trivial for
$d=1$ and follows from   the divergence free property for $d\ge 2$,
we only have to prove {\bf(A0)}, {\bf (A2)} and {\bf (A3)}. Let

$$u=\sum_{k\in\ZZ^d} u_k e_k,\ \ v=\sum_{k\in\ZZ^d} v_k e_k$$ be two
elements in $V_\theta$.

(a) Since $\theta \in [\ff {2\si+1}{ \dd+1}, 1]$ implies $4\si\le 2\theta(\dd+1),$  we have

$$\|u\|_Q^2 =\sum_{k\in\ZZ^d} |u_k|^2|k|^{4\si}\le \ff 1
{\ll_0^\theta} \sum_{k\in\ZZ^d} \ll_0^\theta |u_k|^2
|k|^{2\theta(\dd+1)}=\ff 1 {\ll_0^\theta}\|u\|_{V_\theta}^2.$$ Thus,
{\bf (A0)} holds for $K_1= \ff 1 {\ll_0^\theta}.$

(b) It is easy to see that

\beg{equation}\label{2.2} B(u,v) = \mathbf P
\sum_{l,m\in\ZZ^d, m\ne l}\i
(u_{l-m}\cdot m)v_m  \e_{l}.\end{equation} By H\"older inequality,

\beg{equation*}\beg{split} \|B(u,v)\|_H^2 &\le \sum_{l\in\ZZ^d}\bigg(\sum_{m\in\ZZ^d\setminus\{l\}}|u_{l-m}|\cdot |m|\cdot |v_m|\bigg)^2\\
&\le \sum_{l\in\ZZ^d}\bigg(\sum_{m\in\ZZ^d\setminus\{l\}} |u_{l-m}|^2|m|^{-2\dd}\bigg)\sum_{m\in\ZZ^d} |v_m|^2|m|^{2(\dd+1)}\\
&\le \ff 1 {\ll_0} \bigg(\sum_{m\in\ZZ^d} |m|^{-2\dd}\bigg)\|u\|_H^2\|v\|_V^2.\end{split}\end{equation*}
 Since $\dd>\ff d 2$, we have $\sum_{m\in\ZZ^d} |m|^{-2\dd}<\infty$. Thus,
  {\bf (A2)} holds for some constant $C$.

(c) By (\ref{2.2}), we have

\beg{equation}\label{2.3}\beg{split} &\|B(u,v)\|_Q^2 := \|A^\si B(u,v)\|_H^2\le \sum_{l\in\ZZ^d} |l|^{4\si}\bigg(\sum_{m\in\ZZ^d} |u_{l-m}|\cdot |m|\cdot |v_m|\bigg)^2\\
&\le 2\sum_{l\in\ZZ^d} |l|^{4\si}\bigg(\sum_{|m|>\ff {|l|}2, m\ne l} |u_{l-m}|\cdot |m|\cdot |v_m|\bigg)^2\\
&\qquad
+2\sum_{l\in\ZZ^d} |l|^{4\si}\bigg(\sum_{|m|\le\ff {|l|}2, m\in\ZZ^d} |u_{l-m}|\cdot |m|\cdot |v_m|\bigg)^2:= 2I_1+2I_2.\end{split}\end{equation} By the Schwartz inequality,

$$I_1\le \sum_{l\in\ZZ^d} |l|^{4\si} \bigg(\sum_{|m|>\ff{|l|}2, m\ne l} |u_{l-m}|^2|l-m|^{2(\dd+1)\theta}|m|^{2-2(\dd+1)\theta}\bigg) \sum_{|m|>\ff{|l|}2, m\ne l}|v_m|^2|m|^{2(\dd+1)\theta}|l-m|^{-2(\dd+1)\theta}.$$
Since $\theta\ge \ff{2\si+1}{\dd+1}$ implies that $4\si-2(\dd+1)\theta+2\le 0$,   if $|m|>\ff{|l|}2$ and $|l|\ge 1$ we have

 $$|l|^{4\si}|m|^{-2(\dd+1)\theta+2}\le 4^{(\dd+1)\theta -1}|l|^{4\si-2(\dd+1)\theta +2}\le 4^{(\dd+1)\theta -1}.$$ Therefore,

\beq\label{!}\beg{split} I_1&\le \ff 1 {\ll_0^\theta} 4^{(\dd+1)\theta -1} \|u\|_{V_\theta}^2 \sum_{l\in\ZZ^d}\sum_{|m|>\ff{|l|}2, m\ne l}|v_m|^2|m|^{2(\dd+1)\theta}|l-m|^{-2(\dd+1)\theta}\\
&\le \ff 1 {\ll_0^{2\theta}} 4^{(\dd+1)\theta -1} \bigg(\sum_{m\in\ZZ^d} |m|^{-2(\dd+1)\theta}\bigg)\|u\|_{V_\theta}^2\|v\|_{V_\theta}^2.\end{split}\end{equation}
Similarly,  when $|m|\le \ff {|l|}2 $ we have $|l-m|\ge \ff {|l|}2$ and thus, due to $4\si -2(\dd+1)\theta\le 0$,

$$|l|^{4\si}|l-m|^{-2(\dd+1)\theta}\le 4^{(\dd+1)\theta} |l|^{4\si-2(\dd+1)\theta}\le 4^{(\dd+1)\theta}|m|^{4\si-2(\dd+1)\theta}.$$ Therefore,

\beg{equation*}\beg{split} I_2&\le \sum_{l\in\ZZ^d} |l|^{4\si} \bigg(\sum_{1\le |m|\le \ff {|l|}2} |u_{l-m}|^2|l-m|^{2(\dd+1)\theta}|m|^{2-2(\dd+1)\theta}\bigg)\sum_{1\le |m|\le \ff {|l|}2} |v_m|^2|m|^{2(\dd+1)\theta}|l-m|^{-2(\dd+1)\theta}\\
&\le \ff{4^{(\dd+1)\theta}}{\ll_0^{2\theta}} \bigg(\sum_{m\in\ZZ^d} |m|^{4\si-4(\dd+1)\theta+2}\bigg) \|u\|_{V_\theta}^2\|v\|_{V_\theta}^2\le
\ff{4^{(\dd+1)\theta}}{\ll_0^{2\theta}} \bigg(\sum_{m\in\ZZ^d} |m|^{-2(\dd+1)\theta}\bigg) \|u\|_{V_\theta}^2\|v\|_{V_\theta}^2,\end{split}\end{equation*} where the last step is due to $4\si-2(\dd+1)\theta+2\le 0$ mentioned above. Combining this with (\ref{2.3}) and (\ref{!}), we prove
{\bf (A3)} for the desired $K_2$ which is finite since $\theta\ge \ff {2\si+1}{ \dd+1}$ and $\si>\ff d 4$ imply that
$2(\dd+1)\theta\ge 4\si +1>d.$
\end{proof}

\section{Proofs of Theorem \ref{T1.1} and  consequences}

We first present an exponential estimate of the solution, which will
be used in the proof of Theorem \ref{T1.1}.

\beg{lem} \label{L3.1} In the situation of Theorem \ref{T1.1}, we
have

$$\E\exp\bigg[\ff{\ll_0^2}{2\|Q\|^2}\int_0^t\|X_s^x\|_V^2\d s\bigg]\le
\exp\bigg[\ff{\ll_0^2}{2\|Q\|^2}(\|x\|_H^2+\|Q\|_{HS}^2 t)\bigg],\ \
x\in H, t\ge 0.$$ Moreover, for any $t>0$ and $x\in H$,
$$\E\exp\bigg[\ff{2}{\|Q\|^2\e t}\int_0^t\|X_s^x\|_V^2\d s\bigg]\le
\exp\bigg[\ff{2}{\|Q\|^2t}(\|x\|_H^2+\|Q\|_{HS}^2 t)\bigg].$$
\end{lem}

\beg{proof} (a) Since $\<B(u,v),v\>=0$, by the It\^o formula we have

\beq\label{3.0} \d \|X_t^x\|_H^2 \le -2  \|X_t^x\|_V^2\d t
+\|Q\|_{HS}^2\d t+2\<X_t^x, Q\d W_t\>.\end{equation} Let

$$\tau_n:=\inf\{t\ge 0: \|X_t^x\|_H\ge n\}.$$ By Theorem \ref{T4.1}
below we have $\tau_n\to\infty$ as $n\to\infty$.
 So, for any $\ll>0$ and $n\ge 1$,

\beg{equation*}\beg{split} &\E\,\exp\bigg[\ll \int_0^{t
\land\tau_n}\|X_s^x\|_V^2\d s\bigg]\le \E\exp\bigg[\ff\ll 2
(\|x\|_H^2 +\|Q\|_{HS}^2 t)+ \ll
\int_0^{t\land\tau_n} \<X_s^x, Q\d W_s\>\bigg]\\
&\le \exp\bigg[ \ff\ll 2 (\|x\|_H^2 +\|Q\|_{HS}^2
t)\bigg]\bigg(\E\exp\bigg[ 2\ll^2 \|Q\|^2 \int_0^{t\land\tau_n}
\|X_s^x\|_H^2\d
s\bigg]\bigg)^{1/2}<\infty.\end{split}\end{equation*} Since
$\|\cdot\|_H^2\le \ff 1 {\ll_0}\|\cdot\|_V^2,$ this implies that

 $$\E\,\exp\bigg[\ll \int_0^{t
\land\tau_n}\|X_s^x\|_V^2\d s\bigg]\le\e^{\ff\ll 2 (\|x\|_H^2
+\|Q\|_{HS}^2 t)}\bigg(\E\exp\bigg[\ff{ 2\ll^2 \|Q\|^2}{\ll_0}
\int_0^{t\land\tau_n} \|X_s^x\|_V^2\d s\bigg]\bigg)^{1/2}.$$ Letting
$\ll= \ff{\ll_0^2}{2\|Q\|^2}$, we obtain

$$\E\,\exp\bigg[\ff{\ll_0^2}{2\|Q\|^2} \int_0^{t
\land\tau_n}\|X_s^x\|_V^2\d
s\bigg]\le\exp\bigg[\ff{\ll_0^2}{2\|Q\|^2}(\|x\|_H^2+\|Q\|_{HS}^2
t)\bigg].$$ This proves the first inequality by letting
$n\to\infty.$

(b) Next, due to the first inequality and the Jensen inequality, we
only have to prove the second one for $t\le\ll_0^{-2}.$ In this
case, let
$$\bb(s)= \e^{(\ll_0^2-t^{-1})s},\ \ \ s\in [0,t].$$ By the It\^o formula,  we have

$$\d \|X_s^x\|_H^2\bb(s)=\big\{-2\|X_s^x\|_V^2\bb(s) +\bb'(s)
\|X_s^x\|_H^2+\bb(s)\|Q\|_{HS}^2\big\}\d s +2 \bb(s)\<X_s^x,Q\d
W_s\>.$$ Thus, for any $\ll>0$,

\beg{equation}\label{WW*}\beg{split} &\E\exp\bigg[2\ll\int_0^{t\land
\tau_n}\|X_s^x\|_V^2\bb(s)\d s -\ll \|x\|_H^2 -\ll\|Q||_{HS}^2t\bigg]\\
&\le \E\exp\bigg[2\ll\int_0^{t\land\tau_n} \bb(s) \<X_s^x,Q\d
W_s\>+\ll \int_0^{t\land\tau_n} \bb'(s) \|X_s^x\|_H^2\d s\bigg]\\
&\le \bigg(\E\exp\bigg[2\ll
\int_0^{t\land\tau_n}\|X_s^x\|_V^2\bb(s)\d s\bigg]\bigg)^{1/2}
\bigg(\E\exp\bigg[4\ll \int_0^{t\land \tau_n}\bb(s)\<X_s^x,Q\d
W_s\>\\
&\qquad\qquad \qquad\qquad\qquad\qquad\qquad -2\ll\int_0^{t\land
\tau_n}\|X_s^x\|_H^2\big(\ll_0^2\bb(s)-\bb'(s)\big)\d
s\bigg]\bigg)^{1/2}.\end{split}\end{equation} Note that the first
inequality in the above display implies that

$$\E\exp\bigg[2\ll \int_0^{t\land\tau_n}\|X_s^x\|_V^2\bb(s)\d
s\bigg]<\infty,\ \ n\ge 1.$$ Let $$\ll= \ff 1 {t\|Q\|^2}.$$ By our
choice of $\bb(s)$ and noting that $t\le\ll_0^{-2}$ so that
$\bb(s)\le 1$, we have

$$\ff 1 2 (4\ll)^2\bb(s)^2\|Q\|^2\le
2\ll^2\bb(s)\|Q\|^2\le
2\ll\big(\ll_0^2\bb(s)-\bb'(s)\big).$$Therefore,

$$\E\exp\bigg[4\ll \int_0^{t\land \tau_n}\bb(s)\<X_s^x,Q\d
W_s\>-2\ll\int_0^{t\land
\tau_n}\|X_s^x\|_H^2\big(\ll_0^2\bb(s)-\bb'(s)\big)\d s\bigg]\le
1.$$ Combining this with (\ref{WW*}) for $\ll= (t\|Q\|^2)^{-1},$ we
obtain

$$\E\exp\bigg[\ff{2}{\|Q\|^2\e t}\int_0^{t\land \tau_n}\|X_s^x\|_V^2\d s\bigg]\le
\exp\bigg[\ff{2}{\|Q\|^2t}(\|x\|_H^2+\|Q\|_{HS}^2 t)\bigg].$$ This
completes the proof by letting $n\to\infty$.
\end{proof}

\ \newline \emph{Proof of Theorem \ref{T1.1}.} Simply denote
$X_s= X_s^x$, which solves (\ref{NS}) for $X_0=x$. For given $h\in
V_\theta$ and $\vv>0$, by Theorem \ref{T4.1} below the equation

\beq\label{3.1} \d Y_s= Q\d W_s -\Big\{LY_s +B(X_s) +\ff \vv t
\e^{-Ls}h\Big\}\d s,\ \ Y_0=x+\vv h \end{equation}  has a unique
solution. So,

$$\d (X_s-Y_s) = - L(X_s-Y_s)\d s +\ff \vv t \e^{-Ls} h\d s.$$ This
implies that

\beq\label{3.2} \beg{split} X_s-Y_s&= \e^{-Ls} (X_0-Y_0) +\ff \vv t
\int_0^s \e^{-L(s-r)}\e^{-Lr} h\d r\\
&= \ff{\vv(t-s)}t \e^{-Ls} h =:Z_s,\ \ \ s\in
[0,t].\end{split}\end{equation} Let

$$\eta_s= B(X_s+Z_s)-B(X_s)-\ff \vv t \e^{-Ls}h,$$ which is
well-defined since according to Lemma\ref{L3.1}, $X\in V$ holds $\P\times\d
s$-a.e. Then, by (\ref{3.2}) the equation (\ref{3.1}) reduces to

\beq\label{3.3} \d Y_s= Q\d W_s -\{LY_s +B(Y_s)\}\d s +\eta_s\d
s=Q\d\tt W_s-\{LY_s+B(Y_s)\}\d s,\end{equation} where

$$\tt W_s:= W_s+\int_0^s Q^{-1} \eta_r\d r,\ \ s\in [0,t].$$ By {\bf
(A0)} and  {\bf (A3)}   we have

\beq\label{3.4}\beg{split} \|Q^{-1}\eta_s\|_H^2&\le
\ff{2\vv^2K_1^2}{t^2}\|h\|_{V_\theta}^2 + 2 \|\tt B(X_s,
Z_s)+B(z_s,z_s)\|_Q^2 \\
&\le \vv^2C(t) \big(\|h\|_{V_\theta}^2+\vv^2 \|h\|_{V_\theta}^4
+\|h\|_{V_\theta}^2\|X_s\|_{V_\theta}^2\big).\end{split}\end{equation}
Since $\theta\le 1$ so that $\|\cdot\|_{V_\theta}\le c\|\cdot\|_V$
holds for some constant $c>0$, combining (\ref{3.4}) with Lemma
\ref{L3.1} we concluded that

$$\E \e^{\int_0^t \|\eta_s\|_Q^2\d s}<\infty$$ holds for small
enough $\vv>0.$ By the Girsanov theorem, in this case

$$R_s:=\exp\bigg[-\int_0^s\<Q^{-1}\eta_r,\d W_r\>-\ff 1 2 \int_0^s
\|\eta_r\|_Q^2\d r\bigg],\ \ s\in [0,t]$$ is a martingale and $\{\tt
W_s\}_{s\in [0,t]}$ is the cylindrical Brrownian motion on $H$ under
the probability measure $\R_t\P.$ Combining this with (\ref{3.3})
and the fact that $Y_t=X_t$ due to (\ref{3.2}), for small $\vv>0$ we
have

$$P_t f(x+\vv h)=\E[R_t f(Y_t)] =\E[R_t f(X_t)].$$ Therefore, by the
dominated convergence theorem due to Lemma \ref{L3.1} and
(\ref{3.4}), we conclude that

\beg{equation*}\beg{split} &D_hP_tf(x):=\lim_{\vv\to 0}
\ff{P_tf(x+\vv h)-P_t f(x)}\vv\\
&  = \lim_{\vv\to 0} \E\Big[ \ff{R_t-1}\vv f(X_t)\Big]=
-\E\bigg\{f(X_t)\lim_{\vv\to 0} \int_0^t \Big\<
Q^{-1}\ff{\eta_s}\vv, \d W_s\Big\>\bigg\}\\
&= -\E\bigg\{f(X_t) \int_0^t \Big\<Q^{-1} \Big(\ff{t-s}t \tt
B(\e^{-Ls} h, X_s)-\ff 1 t \e^{-Ls}h\Big), \d
W_s\Big\>\bigg\},\end{split}\end{equation*} where the last step is
due to the bilinear property of $B$, which implies that

\beg{equation*}\beg{split} \ff{\eta_s}\vv &=\ff 1 \vv \tt B(X_s,z_s)+ \ff 1 \vv B(Z_\vv)-\ff 1 t \e^{-Ls}h\\
&=\ff{t-s} t \tt B(X_s, \e^{-Ls}h) -\ff 1 t \e^{-Ls}h+ \ff{\vv(t-s)}t B(\e^{-Ls}h,\e^{-Ls}h).\end{split}\end{equation*}\qed

\ \newline\emph{Proof of Corollary \ref{C1.2}.}  (1) By (\ref{B}) and
the Schwartz inequality, for any $h$ with $\|h\|_{V_\theta}\le 1,$
we have

\beq\label{3.5}\beg{split} |D_h P_t f(x)|^2 &\le
(P_tf(x))^2\E\int_0^t \Big\|\ff 1 t \e^{-Ls}h-\ff{t-s}t \tt B(X_s^x,
h)\Big\|_Q^2\d s\\
&\le 2 (P_tf^2(x)) \bigg\{\ff {K_1} t  + \E \int_0^t \|\tt B(X_s^x,
h)\|_Q^2\d s\bigg\},\end{split}\end{equation} where the last step is
due to the fact that  {\bf (A0)} implies

\beq\label{3.*} \|\e^{-Ls}h\|_Q^2\le K_1
\|\e^{-Ls}h\|_{V_\theta}^2\le K_1\|h\|_{V_\theta}^2.\end{equation}
Next, by {\bf (A3)} and $\theta\le 1$ we have

\beq\label{3.**} \|\tt B(X_s^x,h)\|_Q^2\le 4K_2
\|h\|_{V_\theta}^2\|X_s^x\|_{V_\theta}^2\le
\ff{4K_2}{\ll_0^{1-\theta}} \|X_s^x\|_V^2.\end{equation} Combining
this with (\ref{3.0}) we obtain

$$\E\int_0^t \|\tt B(X_s^x,h)\|_Q^2\d s \le \ff{2
K_2}{\ll_0^{1-\theta}} \big(\|x\|_H^2+\|Q\|_{HS}^2 t\big).$$ The
proof of (1) is completed by this and (\ref{3.5}).

(2) Let $f\ge 0$ and $h$ be such that $\|h\|_{V_\theta}\le 1.$ Let

$$M_t=  \int_0^t \Big\<Q^{-1} \Big(\ff{t-s}t \tt
B(\e^{-Ls} h, X_s)-\ff 1 t \e^{-Ls}h\Big), \d W_s\Big\>.$$ By
(\ref{B}) and the Young inequality (see e.g. \cite[Lemma
2.4]{ATW09}),

\beq\label{3.6} |D_h P_tf(x)|\le \dd\big\{P_t(f\log f)-(P_t f)\log
P_t f\big\}(x) +\big\{\dd  \log\E \e^{\ff 1 \dd M_t}\big\}P_tf(x) ,\
\ \dd>0.\end{equation} Since by (\ref{3.*}) and (\ref{3.**}) we have

\beg{equation}\label{WW2}\beg{split} \<M\>_t &= \int_0^t \Big\|\ff 1 t
\e^{-Ls} h-\ff{t-s}t \tt B(X_s^x, h)\Big\|_Q^2 \d s \\
&\le \ff{2K_1}t + \ff{4 K_2}{\ll_0^{1-\theta}} \int_0^t
\|X_s^x\|_V^2\d s,\end{split}\end{equation} it follows from
Lemma \ref{L3.1} that for any $\dd\ge \dd_0:=
4\ss{K_2}\,\|Q\|\ll_0^{(\theta-3)/2},$

\beg{equation*}\beg{split} \E \exp\bigg[\ff 1 \dd M_t\bigg]&\le
\bigg(\E \exp\bigg[\ff 2 {\dd^2} \<M\>_t\bigg]\bigg)^{1/2} \le
\bigg(\E\exp\bigg[\ff 2
{\dd_0^2}\<M\>_t\bigg]\bigg)^{\dd_0^2/(2\dd^2)}\\
&\le \exp\bigg[\ff{2K_1}{\dd^2t}\bigg]
\bigg(\E\exp\bigg[\ff{8K_2}{\dd_0^2\ll_0^{1-\theta}}\int_0^t
\|X_s^x\|_V^2\d s\bigg]\bigg)^{\dd_0^2/(2\dd^2)}\\
&=\exp\bigg[ \ff{2K_1}{\dd^2t}\bigg]
\bigg(\E\exp\bigg[\ff{\ll_0^2}{2\|Q\|^2}\int_0^t \|X_s^x\|_V^2\d
s\bigg]\bigg)^{\dd_0^2/(2\dd^2)}\\
&\le
\exp\bigg\{\ff{2K_1}{\dd^2t}+\ff{\ll_0^2\dd_0^2}{4\dd^2\|Q\|^2}(\|x\|_H^2+\|Q\|_{HS}^2
t)\bigg\}\\
&=\exp\bigg\{\ff 2 {\dd^2}\Big(\ff{K_1}t
+\ff{2K_2}{\ll_0^{1-\theta}}(\|x\|_H^2+\|Q\|_{HS}^2t)\Big)\bigg\}.\end{split}\end{equation*}
Combining this with (\ref{3.6}) we prove (2).

(3) According to e.g. \cite[proof of Proposition 4.1]{DRW09}), the
$V_\theta$-strong Feller property of $P_t$ follows from the claimed
Harnack inequality, which we prove below   by using an argument in
\cite[Proof of Theorem 1.2]{ATW09}. Let $x\ne y$ be such that

\beq\label{?}\|x-y\|_{V_\theta}\le \ff{\aa-1}{\aa\dd_0}\ \text{for}\  \dd_0:=
\ff{4\|Q\|\ss {K_2}}{\ll_0^{(3-\theta)/2}}.\end{equation} Let

$$\bb_s=1+s(\aa-1),\ \ \gg_s= x+s(y-x),\ \ s\in [0,1].$$ We have

\beg{equation*}\beg{split} &\ff{\d}{\d s} \log (P_t
f^{\bb(s)})^{\aa/\bb(s)}(\gg_s)\\
&= \ff{\aa(\aa-1)}{\bb(s)^2}\cdot \ff{P_t(f^{\bb(s)}\log f^{\bb(s)})
-(P_tf^{\bb(s)})\log P_t f^{\bb(s)}}{P_t f^{\bb(s)}}(\gg_s)+\ff{\aa
D_{y-x} P_t f^{\bb(s)}}{\bb(s)P_t f^{\bb(s)}}(\gg_s)\\
&\ge \ff{\aa\|x-y\|_{V_\theta}}{\bb(s)P_t f^{\bb(s)}(\gg_s)}
\bigg\{\ff{\aa-1}{\bb(s)\|x-y\|_{V_\theta}}\Big(P_t(f^{\bb(s)}\log
f^{\bb(s)}) -(P_tf^{\bb(s)})\log P_t f^{\bb(s)}\Big)(\gg_s)\\
&\qquad\qquad\qquad\qquad\qquad\qquad\qquad\qquad\qquad\qquad\qquad\qquad\qquad-\|DP_t
f^{\bb(s)}(\gg_s)\|_{V_\theta}^*\bigg\}.\end{split}\end{equation*}Therefore,
applying (2) to

$$\dd:= \ff{\aa-1}{\bb(s) \|x-y\|_{V_\theta}}$$ which is larger than
$\dd_0$ according to (\ref{?}), we obtain

\beg{equation*}\beg{split}  &\ff{\d}{\d s} \log (P_t
f^{\bb(s)})^{\aa/\bb(s)}(\gg_s) \ge
-\ff{2\aa\|x-y\|_{V_\theta}}{\dd\bb(s)}\bigg\{\ff{K_1}t
+\ff{2K_2}{\ll_0^{1-\theta}} (\|\gg_s\|_H^2 +\|Q\|_{HS}^2 t)\bigg\}\\
&\ge -\ff{2\aa\|x-y\|^2_{V_\theta}}{\aa-1}\bigg\{\ff{K_1}t
+\ff{2K_2}{\ll_0^{1-\theta}} (\|x\|_H^2\lor\|y\|_H^2 +\|Q\|_{HS}^2
t)\bigg\}.\end{split}\end{equation*}Integrating over $[0,1]$ w.r.t.
$\d s$, we derive the desired Harnack inequality.\qed

\ \newline\emph{Proof of Corollary \ref{C1.3}.}   Since $u\mapsto
\|u\|_V^2$ is a compact function on $H$, i.e. for any $r>0$ the set
$\{u\in H: \|u\|_V\le r\}$ is relatively compact in $H$, (\ref{3.0})
implies the existence of the invariant probability measure
satisfying (1) by a standard  argument (see e.g. \cite[Proof of
Theorem 1.2]{W07}). Moreover, any invariant probability measure
$\mu$ satisfies $\mu(\|\cdot\|_V^2)<\infty$, hence, $\mu(V)=1.$
Below, we assume $\theta\in (0,1)$ and prove (1) and (2)
repsectively.

(1) Let $\mu$ be an invariant probability measure, we first prove it
has full support on $\mu$. $$r_0=
\ff{\ll_0^{(3-\theta)/2}}{8\|Q\|\ss{K_2}}.$$ By Corollary
\ref{C1.2}(3) for $\aa=2$, for any fixed $t>0$ there exists a constant
$C(t)>0$ such that

$$(P_t f(x))^2\le (P_t f^2(y))\e^{C(t)(\|x\|_H^2+\|y\|_H^2)},\ \
\|x-y\|_{V_\theta}\le r_0.$$ Applying this inequality $n$ times, we
may find a constant $c(t,n)>0$ such that

\beq\label{HH}(P_t f(x))^{2n}\le (P_t
f^{2n}(y))\e^{C(t,n)(\|x\|_H^2+\|y\|_H^2)},\ \ \|x-y\|_{V_\theta}\le
nr_0.\end{equation} Since $V$ is dense in $V_\theta$, to prove that
$\mu$ has full support on $V_\theta$,   it suffices to show that

\beq\label{S} \mu(B_\theta(x,\vv))>0,\ \ \ x\in V, \vv>0
\end{equation} holds for  $
B_\theta(x,\vv):=\{y: \|y-x\|_{V_\theta}<\vv\}.$ Since
$\mu(V_\theta)=1,$ there exists $n\ge 1$ such that
$\mu(B_\theta(x,nr_0))>0.$ Applying (\ref{HH}) to
$f=1_{B_\theta(x,\vv)}$ we obtain

$$\P(\|X_t^x-x\|_{V_\theta}<\vv)^{2n}\int_{B_\theta(x,nr_0)}\e^{-C(t,n)(\|x\|_H^2+\|y\|_H^2)}\mu(\d
y)\le\mu(B_\theta(x,\vv)).$$ So,  if $\mu(B_\theta (x,\vv))=0$ then

\beq\label{CC}\P(\|X_t^x-x\|_{V_\theta}\ge \vv)=1,\ \ \
t>0.\end{equation} To see that this is impossible, let us observe
that for any $m\ge 1$ there exists a constant $c(m)>0$ such that

\beq\label{!!}\|\cdot\|_{V_\theta}^2\le c(m)\|\cdot\|_H^2+ \ff 1
{(\ll_0m)^{1-\theta}}\|\cdot\|_V^2\end{equation} holds. Moreover,
using $\<\cdot,\cdot\>$ to denote the duality w.r.t $H$,  we have

\beg{equation*}\beg{split} &2\<X_t^x-x,LX_t^x\> = 2\|X_t^x-x\|_V^2 +2 \<X_t^x-x,Lx\> \\
&\ge 2\|X_t^x-x\|_V^2-2\|X_t^x-x\|_V\|x\|_V \ge \|X_t^x-x\|_V^2 -
\|x\|_V^2\end{split}\end{equation*} and due to {\bf (A1)} and {\bf
(A2)},
$$2\<X_t^x-x, B(X_t^x)\>=-2\<x,B(X_t^x)\>\le 2C\|x\|_H\|X_t^x\|_V\|X_t^x\|_H\le
\ff 1 2 \|X_t^x-x\|_V^2 + c_1+c_2\|X_t^x\|_H^2$$ holds for some
constants $c_1,c_2$ depending on $x$. Therefore, by the It\^o
formula for $\|X_t^x-x\|_H^2$, we arrive at

\beg{equation*}\beg{split} \d\|X_t^x-x\|_H^2 &= \big\{\|Q\|_{HS}^2-2\<X_t^x-x,LX_t^x\>+ 2\<X_t^x-x, B(X_t^x)\>\big\}\d t+ 2\<X_t^x-x,Q\d W_t\>\\
&\le -\ff 1 2 \|X_t^x-x\|_V^2\d t + (c_3+c_2\|X_t^x\|_H^2)\d t +2\<X_t^x-x,Q\d W_t\>\end{split}\end{equation*} for some constant
$c_3>0$. Since by Theorem \ref{T4.1} below $\E \sup_{t\in [0,1]}\|X_t\|_H^2<\infty$, this and the continuity of $X_s^x$ in $s$ imply
 
 $$\lim_{t\to 0} \ff 1 t \int_0^t \|X_s^x-x\|_H^2\d s=0$$ and

$$\E\int_0^t \|X_t^x-x\|_V^2\d s\le c_0t,\ \ \ t\in[0,1]$$  for
some constant $c_0>0$.    Combining these with (\ref{!!}),  we conclude
 that

$$\limsup_{t\to 0}\ff 1 t \int_0^t \E\|X_s^x-x\|_{V_\theta}^2 \d s
\le \ff {c_0} {(\ll_0m)^{1-\theta}},\ \ \ m\ge 1.$$ Letting
$m\to\infty$ we obtain

$$\lim_{t\to 0}\ff 1 t \int_0^t \E\|X_s^x-x\|_{V_\theta}^2 \d s=0.$$
this is contractive to (\ref{CC}).

Next, if the invariant probability measure is not unique, we may
take two different extreme elements $\mu_1,\mu_2$ of the set of all
invariant probability measures. It is well-known that $\mu_1$ and
$\mu_2$ are singular with each other. Let $D$ be a $\mu_1$-null set,
since $\mu_1$ has full support on $V_\theta$ and $P_t 1_D$ is
continuous and $\mu_1(P_t 1_D)= \mu_1(D)=0$, we have $P_t1_D\equiv
0$. Thus, $\mu_2(D)=\mu_2(P_t1_D)=0$. This means that $\mu_2$ has to
be absolutely continuous w.r.t. $\mu_1$, which is contradictive to
the singularity of $\mu_1$ and $\mu_2$.

(2)  As observe above that   $P_t 1_D\equiv 0$ for any $\mu$-null
set $D$. So, $P_t$ has a transition density $p_t(x,y)$ w.r.t. $\mu$
on $V_\theta$.  Next, let $f\ge 0$ such that $\mu(f^\aa)\le 1$. By
the Harnack inequality in Corollary \ref{C1.2}(3), we have

$$(P_t f(x))^\aa \int_{B_\theta(x,r_0)}
\exp\bigg[-\ff{2\aa\|x-y\|_{V_\theta}^2}{\aa-1}\Big\{\ff{K_1}t
+\ff{2K_2}{\ll_0^{1-\theta}} \big(\|x\|_H^2\lor\|y\|_H^2+
\|Q\|_{HS}^2t\big)\Big\}\bigg]\mu(\d y)\le 1.$$ Then the desired
estimate on $\int p_t(x,z)^{(\aa+1)/\aa a}\mu(\d z)$ follows by
taking

$$f(\cdot)= p_t(x,\cdot).$$

\ \newline \emph{Proof of Theorem \ref{T1.4}.} (1) Let $M_t$ be in
the proof of Corollary \ref{C1.2} (2). By (\ref{WW2}), for $\dd>0$
we have

\beg{equation*}\beg{split} &\E\exp\Big[\ff {M_t}\dd\Big]\le \Big(\E
\exp\Big[\ff{2\<M\>_t}{\dd^2}\Big]\Big)^{1/2}\\
& \le \exp\Big[\ff {2K_1}{\dd^2 t}\Big]
\bigg(\exp\bigg[\ff{8K_2}{\ll_0^{1-\theta}\dd^2}\int_0^t
\|X_s^x\|_V^2\d s\bigg]\bigg)^{1/2}.\end{split}\end{equation*} If
$t\le t_\dd$ then

$$\ff{8K_2}{\ll_0^{1-\theta}\dd^2}\le \ff{2\ll_0^2}{\|Q\|^2\e t},$$
so that by the Jensen inequality and the second inequality in Lemma
\ref{L3.1},

\beg{equation*}\beg{split}\E\exp\Big[\ff
{M_t}\dd\Big]&\le\exp\Big[\ff {2K_1}{\dd^2 t}\Big]
\bigg(\exp\bigg[\ff{2\ll_0^2}{\|Q\|^2\e t}\int_0^t \|X_s^x\|_V^2\d
s\bigg]\bigg)^{\ff{2K_2\|Q\|^2\e t}{\dd^2\ll_0^{3-\theta}}}\\
&\le \exp\Big[\ff {2K_1}{\dd^2
t}+\ff{4K_2\e}{\dd^2\ll_0^{1-\theta}}\Big],\ \ t\le
t_\dd.\end{split}\end{equation*} Combining this with (\ref{3.6}) we
prove the desired gradient estimate for $t\le t_\dd$. By the
gradient estimate   for $t=t_\dd$ and the semigroup property, when
$t>t_\dd$ we have

\beg{equation*}\beg{split}&\|DP_tf(x)\|_{V_\theta^*}=\|DP_{t_\dd}(P_{t-t_\dd}f)(x)\|_{V_\theta^*}\le
\dd \big\{P_{t_\dd}\big((P_{t-t_\dd}f)\log P_{t-t_\dd}f\big)\\
&\qquad -(P_t f)\log P_t f\big\}(x) + \ff 2 \dd
\Big\{\ff{K_1}{t_\dd} +\ff{2K_2\e}{\ll_0^{1-\theta}}\big(\|x\|_H^2
+\|Q\|_{HS}^2t\big)\Big\}P_tf(x).\end{split}\end{equation*} This
implies the desired gradient estimate for $t>t_\dd$ since due to the
Jensen inequality
$$P_{t_\dd}\big((P_{t-t_\dd}f)\log P_{t-t_\dd}f\big)\le P_t f\log
f.$$

(2) Repeating the proof of Corollary \ref{C1.3} (3) using the
inequality in Theorem \ref{T1.4} (1) instead of Corollary \ref{C1.2}
(2) for $\dd= \ff{\aa-1}{\bb(s)\|x-y\|_{V_\theta}}$, we obtain

$$\ff{\d}{\d s} \big(\log P_t f^{\bb(s)}\big)^{\aa/\bb(s)}\ge
-\ff{2\aa\|x-y\|_{V_\theta}^2}{\aa-1}\bigg\{\ff{K_1}{t\land
t_\dd}+\ff{2K_2
\e}{\ll_0^{1-\theta}}\big(\|x\|_H^2\lor\|y\|_H^2+\|Q\|_{HS}^2
t\big)\bigg\}.$$ This completes the proof by integrating over
$[0,1]$ w.r.t. $\d s$ and noting that

$$ t_\dd = \ff{\dd^2\ll_0^{3-\theta}}{4\|Q\|^2\e K_2}\ge
\ff{(\aa-1)^2\ll_0^{3-\theta}}{4\aa^2\|Q\|^2 K_2\e
\|x-y\|_{V_\theta}^2}$$ since

$$\dd= \ff{\aa-1}{\bb(s)\|x-y\|_{V_\theta}}\ge
\ff{\aa-1}{\aa\|x-y\|_{V_\theta}}.$$\qed

\section{Appendix} We aim to verify the existence and uniqueness of
the solution to (\ref{E}) by using the main result of \cite{LR}.

\beg{thm} \label{T4.1} Assume {\bf (A1)} and {\bf (A2)}. For any
$X_0\in H$ the equation $(\ref{E})$ has a unique solution $X_t$,
which is a continuous Markov process on $H$ such that

$$\E\bigg(\sup_{t\in [0,T]}\|X_t\|_H^p+ \int_0^T \|X_t\|_V^2\d
t\bigg)<0$$ holds for any $p>1$ and $\P$-a.s.

$$X_t =X_0-\int_0^t(LX_s +B(X_s))\d s+QW_t,\ \ t\ge 0$$ holds on
$H$. \end{thm}

\beg{proof} Let $V^*$ be the dual space of $V$ w.r.t. $H$. Then for
any $v\in V$,

$$A(v):= -(Lv+B(v))\in V^*.$$ It suffices to verify assumptions
(H1)-(H4) in \cite[Theorem 1.1]{LR} for the functional $A$. The
hemicontinuity assumption H1) follows immediately form the bilinear
property of $B$. Next, by {\bf (A2)} and the bilinear property of
$B$, we have

\beg{equation*}\beg{split} _{V^*}\<A(v_1)-A(v_2),v_1-v_2\>_V
&=-\|v_1-v_2\|_V^2 +\|B(v_2-v_1,v_1),v_1-v_2\>\\
&\le -\|v_1-v_2\|_V^2
+C\|v_1-v_2\|_H^2\|v_1\|_V^2.\end{split}\end{equation*} So, the
assumption (H2) in \cite{LR} holds for $\rr(v):= c\|v\|_V^2.$
Moreover, by {\bf (A1)} we have

$$_{V^*}\<A(v),v\>_V\le -\|v\|_V^2.$$ Thus, the coercivity assumption
(H3) in \cite{LR} holds for $\theta=1, \aa=2, K=0$ and $f=$constant.
Finally, {\bf (A2)} implies that

$$\|A(v)\|_{V^*}^2 \le 2 \|v\|_V^2 + 2 \|L^{-1/2}B(v)\|_H^2\le 2
\|v\|_V^2 +\ff {2c} {\ll_0} \|v\|_H^2\|v\|_V^2.$$ Therefore, the
growth condition (H4) in \cite{LR} holds for some constant $f,K>0$
and $\aa=\bb=2.$\end{proof}
 \beg{thebibliography}{99}

\bibitem{ATW06}  M. Arnaudon, A. Thalmaier, F.-Y. Wang,
  \emph{Harnack inequality and heat kernel estimates
  on manifolds with curvature unbounded below,} Bull. Sci. Math. 130(2006), 223--233.

\bibitem{ATW09} M. Arnaudon, A. Thalmaier and F.-Y. Wang,
  \emph{Gradient estimates and Harnack inequalities on non-compact Riemannian manifolds,}
    Stoch. Proc. Appl. 119(2009), 3653--3670.

\bibitem{DD} G. Da Prato, A. Debussche,  \emph{Ergodicity for the 3D stochastic Navier-Stokes equations,}
J. Math. Pures Appl. 82(2003), 877--947.

\bibitem{DRW09} G. Da Prato, M. R\"ockner and F.-Y. Wang,
  \emph{Singular stochastic equations on Hilbert spaces:  Harnack inequalities for their
  transition semigroups,}  J. Funct. Anal. 257(2009), 992--1017.

\bibitem{FR} F. Flandoli, M. Romito, \emph{Markov selections for the 3D stochastic Navier-Stokes equations,}
Probab. Theory Relat. Fields 140(2008), 407--458.

 \bibitem{HM} M. Hairer, J. C. Mattingly, \emph{Ergodicity of the 2D Navier-Stokes equations with degenerate
  stochastic forcing,} Ann. Math. 164(2006), 993--1032.

  \bibitem{L} W. Liu, Doctor-Thesis, Bielefeld University, 2009.

\bibitem{LR} W. Liu, M. R\"ockner, \emph{SPDE in Hilbert space
  with locally monotone coefficients,} arXiv:1005.0632v1, 2010.

\bibitem{LW08} W. Liu and F.-Y. Wang, \emph{Harnack inequality and strong Feller
  property for stochastic fast diffusion equations,} J. Math. Anal. Appl.
  342(2008), 651--662.

\bibitem{OY} S.-X. Ouyang, Doctor-Thesis, Bielefeld University, 2009.

 \bibitem{ORW} S.-X. Ouyang, M. R\"ockner and F.-Y. Wang,
 \emph{Harnack inequalities and applications for Ornstein-Uhlenbeck semigroups with
 jump,} arXiv:0908.2889

\bibitem{RX} M. Romito, L. Xu, \emph{HErgodicity of the 3D stochastic Navier-Stokes equations driven by mildly degenerate
noise,} 2009 Preprint.

\bibitem{T} R. Temam, \emph{Navier-Stokes Equations and Nonlinear Functional Analysis (2nd Ed),}
CBMS-NSF Regional Conference Series in Appl. Math. V66, SIAM,
Philadelphia, 1995.

\bibitem{W97} F.-Y. Wang, \emph{Logarithmic Sobolev
inequalities on noncompact Riemannian manifolds,} Probab. Theory
Relat. Fields 109(1997), 417--424.

\bibitem{W07}  F.-Y. Wang, \emph{ Harnack inequality and applications for stochastic generalized porous media equations,}
  Ann. Probab. 35(2007), 1333--1350.

 \end{thebibliography}
\end{document}